\RequirePackage{fix-cm}
\documentclass[smallextended]{svjour3}       % onecolumn (second format)
\smartqed  % flush right qed marks, e.g. at end of proof
\usepackage[dvips]{graphicx}
\usepackage[dvips]{color}
\usepackage{amsmath}
\usepackage{amssymb}
\usepackage{epsfig}
\usepackage{mathrsfs}
\usepackage{amsfonts}
\def\pn{\par\smallskip\noindent}
\def\proof{\pn {Proof.} }
\def\eproof{\hfill \quad{$\Box$}\smallskip}

\newtheorem{lem}{Lemma}

\newtheorem{thm}{Theorem}
\newtheorem{cor}{Corollary}
\newtheorem{pro}{Proposition}

\begin{document}

\title{On the p-regularized trust region subproblem
\thanks{This research was supported by Chinese NSF grants 10831006, 11021101, and CAS grant kjcx-yw-s7, by Taiwan NSC 98-2115-M-006-010-MY2, by Beijing Higher Education Young Elite
Teacher Project 29201442, and by the fund of State Key Laboratory of
Software Development Environment SKLSDE-2013ZX-13.}}

%\author{Shu-Cherng Fang\footnotemark[2] \and David Yang Gao\footnotemark[3]
% \and Gang-Xuan Lin\footnotemark[4] \and Ruey-Lin Sheu\footnotemark[4]~\footnotemark[5]\and Wen-Xun Xing\footnotemark[6]}

\titlerunning{The $p-$regularized subproblem}        % if too long for running head

\author{Yong Hsia \and Ruey-Lin Sheu\and Ya-xiang Yuan}

%\authorrunning{Y. Hsia et al} % if too long for running head

\institute{Y. Hsia \at State Key Laboratory of Software Development
              Environment, LMIB of the Ministry of Education, School of
Mathematics and System Sciences, Beihang University,  Beijing
100191,  P. R. China
              \email{dearyxia@gmail.com}\\
  R. L. Sheu \at
              Department of Mathematics, National Cheng Kung University, Taiwan
              \email{rsheu@mail.ncku.edu.tw}\\
  Y. Yuan \at
              State Key Laboratory of Scientific/Engineering Computing,
              Institute of Computational Mathematics
and Scientific/Engineering Computing, The Academy of Mathematics and
Systems Sciences, Chinese Academy of Sciences, P.O. Box 2719,
Beijing 100080, P. R. China
              \email{yyx@lsec.cc.ac.cn}            }

\date{Received: date / Accepted: date}

\maketitle

\begin{abstract}

The $p$-regularized subproblem (p-RS) is a regularisation technique
in computing a Newton-like step for unconstrained optimization,
which globally minimizes a local quadratic approximation of the
objective function while incorporating with a weighted
regularisation term $\frac{\sigma}{p} \|x\|^p$. The global solution
of the $p$-regularized subproblem for $p=3$, also known as the cubic
regularization, has been characterized in literature. In this paper,
we resolve both the global and the local non-global minimizers of
(p-RS) for $p>2$ with necessary and sufficient optimality
conditions. Moreover, we prove a parallel result of Mart\'{\i}nez
\cite{Mar} that the (p-RS) for $p>2$, analogous to the trust region
subproblem, can have at most one local non-global minimizer. When
the (p-RS) is subject to a fixed number $m$ additional linear
inequality constraints, we show that the uniqueness of the local
solution of the (p-RS) (if exists at all), especially for $p=4$, can
be applied to solve such an extension in polynomial
time.%Our theoretic results set a foundation for
%developing a whole class of global convergence $p$-regularised
%Newton's algorithm for any $p>2$, and provide an alternative to
%problems wherever the trust region method may apply.

\end{abstract}

\section{Introduction}

For an unconstrained optimization problem to minimize $f$ over $\Bbb
R^n$, Newton's method has an attractive local convergence property
near a second order critical point. Ensuring the global convergence
for Newton's method with an analyzable computational complexity,
however, requires modifications to guarantee a {\it sufficient}
descent at each step. Unlike the Levenberg-Marquardt type of methods
or most quasi-Newton methods which always maintain a
positive-definite approximate Hessian of $f$, the $p$-regularized
subproblem minimizes globally the second order Taylor's polynomial
of $f$ plus a weighted (by $\sigma$) higher order regularization
term. The subproblem takes the following model
\begin{equation}
({\rm p-RS})~~\min_{x\in \Bbb R^n} \left\{
g(x)= \frac{1}{2}x^THx+c^Tx+\frac{\sigma}{p} \|x\|^p \right\}, \label{d-well}
\end{equation}
where $\sigma>0$, $p>2$, and $H$ is the Hessian of $f$ at any
iterate, regardless of its definiteness. It is often assumed that
$f$ is smooth enough to have a symmetric Hessian and to obtain the
desire global convergence. At each iterate, if the global minimizer
of ({\rm p-RS}) renders a satisfactory decrease in the value of $f$,
it is accepted; but rejected otherwise with an increase in $\sigma$
to enhance the regularization force.

In literature, ({\rm p-RS}) with $p = 3$ is known as the cubic
regularization which is the most common choice among all others. The
idea of the cubic regularization was first due to Griewank \cite{Gr}
and later was considered by many authors with thorough global
convergence and complexity analysis. See Nesterov and Polyak
\cite{Ne}; Weiser Deuflhard and Erdmann \cite{W}; and Cartis, Gould
and Toint \cite{Cartis}. When $p=4$, (p-RS) reduces to a form of the
double well potential function which has many applications in solid
mechanics and quantum mechanics \cite{Fang,Xia}. Gould, Robinson and
Thorne \cite{G10} studied ({\rm p-RS}) for a general $p>2$ in
comparison with the the trust-region subproblem
\begin{eqnarray}
({\rm TRS})~~&\min & \frac{1}{2}x^THx+c^Tx \label{tr:1}\\
&{\rm s.t.}& \Vert x\Vert^2 \le \Delta, ~x \in\Bbb R^n. \label{tr:2}
\end{eqnarray}

Our paper characterizes (p-RS) completely for any $p>2$ by extending
(i) the necessary and sufficient global optimality conditions for
$p=3$ in \cite{Cartis}; (ii) the analysis using the secular function
(to be specified later) for $p=4$ in \cite{Xia}; and (iii) a
necessary global optimality condition for $p>2$ in \cite{G10}. Some
generalization is, nevertheless, non-trivial in mathematical skills.
We summarize the main results as follows.

\begin{itemize}
\item[$\bullet$] Theorem 1 of the paper ({\it cf.} Theorem 3.1 in \cite{Cartis} for $p=3$;
Theorem 2 in \cite{G10} for the necessary part of $p>2$): The
point $x^*$ is a global minimizer of (p-RS) for $p
> 2$ if and only if
\begin{equation*}
(H+\sigma\|x^*\|^{p-2}I)x^*=-c\ ;~~~H+\sigma\|x^*\|^{p-2}I\succeq 0.  %\label{R:1}
\end{equation*}
%with $H+\lambda^*I\succeq 0$ and $\lambda^*=\sigma
%\|x^*\|^{p-2}$.

\item[$\bullet$] Theorem 2 ({\it cf.} the trust region subproblem in \cite{HS}):
Let $k$ be the multiplicity of the smallest eigenvalue
$\alpha_1$ of $H$, i.e.,
\[
\alpha_1=\ldots=\alpha_k<\alpha_{k+1}\le\ldots\le\alpha_n.
\]
Then, the set of the global minimizers of (p-RS) is either a
singleton or a $k$-dimensional sphere centered
    at
    $(0,\cdots,0,-\frac{c_{k+1}}{\alpha_{k+1}-\alpha_1},\cdots,-\frac{c_{n}}{\sigma_{n}-\sigma_1})$
    with the radius
    $\sqrt{\left(\frac{\alpha_1}{\sigma}\right)^{\frac{2}{p-2}}-\sum_{i=k+1}^n\frac{c_i^2}{(\alpha_i-\alpha_1)^2}}$.
\item[$\bullet$]
Theorem 3 ({\it cf.} Theorem 2 in \cite{Xia} for $p=4$): The
point $\underline{x}$ is a local-nonglobal minimizer of ({\rm
p-RS}) for $p>2$ if and only if
\begin{equation}
\underline{x}
=-
\left( H+\sigma\underline{t}^* I\right)^{-1}c,\label{tw}
\end{equation}
%\underline{x}
%=\left(\frac{-c_1}{\sigma\underline{t}^*+\alpha_1},\ldots,\frac{-c_n}
%{\sigma\underline{t}^* +\alpha_n} \right),
where $\underline{t}^*$ is a root of the secular function
\begin{equation}
 h(t)=\|\left(H+\sigma t  I\right)^{-1}c\|^2-t^{\frac{2}{p-2}},
 ~~t\in\left(\max\left\{- \frac{\alpha_2}{\sigma},0\right\}, - \frac{\alpha_1}
 {\sigma}\right)\label{varphi00}
\end{equation}
such that $h'(\underline{t}^*)>0$.
\item[$\bullet$] Theorem 4 ({\it cf.} the trust region subproblem in \cite{Mar}; the double
well potential function in \cite{Xia}): The subproblem (p-RS)
with $p>2$ has at most one local non-global minimizer.
\end{itemize}
Notice that, the secular function for (TRS) ({\it cf.} $h(t)$ in
(\ref{varphi00})) is defined by $$ \phi(\lambda)=\|(H+\lambda
I)^{-1}c\|^2.$$ Mart\'{\i}nez \cite{Mar} proved that, if
$\underline{x}$ is a local-nonglobal minimizer of (TRS), then
$\underline{x}$ satisfies $(H+\lambda^*I)\underline{x}=-c$ with
$\lambda^*\in (-\alpha_2,-\alpha_1)$, $\lambda^*\ge 0$ and
$\phi'(\lambda^*)\ge 0$. However, to the best of our knowledge, the
necessary condition $\phi'(\lambda^*)\ge 0$ is not known to be
sufficient for (TRS) or not.

Finally, as an application, we study (p-RS) subject to $m$ linear
inequality constraints of the following form:
\begin{eqnarray}
({\rm p-RS}_m)~&\min & \frac{1}{2}x^THx+c^Tx+\frac{\sigma}{p} \|x\|^p  \label{RSm:1}\\
&{\rm s.t.}&  l_i\le a_i^Tx\le u_i,~i=1,\ldots,m, \label{RSm:2}
\end{eqnarray}
where $l_i\le u_i\in \Bbb R$ for $i=1,\ldots,m$. We first show that
the NP-hard $k$-dispersion-sum problem
\begin{eqnarray}
{\rm(KDSP)}~~d^*=&\min &x^T(-D)x \label{p2:1}\\
&{\rm s.t.}&e^Tx=k,~x\in \{0,1\}^n\label{p2:2}
\end{eqnarray}
can be reduced to a special case of ({\rm p-RS}$_{n+1})$ with $p=4$.
It indicates that solving the class of subproblems
$\bigcup_{m>n}$({\rm p-RS}$_{m})$ with $p=4$ is also NP-hard.
However, for any fixed $m$, by Theorem 4 that there is at most one
local non-global minimizer for ({\rm p-RS}$_{m})$ with $p=4$, we
show that it can be solved in polynomial time. Notice that there is
an analogy called the extended trust region subproblem which adds
linear inequality constraints to (TRS). Polynomial solvability has
been recently proved by Bienstock and Michalka \cite{Bienstock}, and
independently by Hsia and Sheu \cite{HS}.

{\bf Notations.} Let $v(\cdot)$ denote the optimal value of problem
$(\cdot)$. For any symmetric matrix $P\in\Bbb R^{n\times n}$,
$P\succ (\succeq) 0$ means that $P$ is positive (semi)definite. The
determinant of $P$ is denoted by $\det(P)$ whereas the identity
matrix of order $n$ by $I$. For a vector $x\in\Bbb R^n$, Diag$(x)$
is a diagonal matrix with diagonal components being
$x_1,\ldots,x_n$. For a number $\beta \in\Bbb R$,
sign$(\beta)=\frac{\beta}{|\beta|}$ if $\beta \neq 0$, otherwise
sign$(\beta)=0$.

%The global solutions of (TRS) were fully characterized in the early 1980s by Gay \cite{Ga}, Sorensen \cite{So} and Mor\'{e} and Sorensen \cite{Mo}:
% \\[+5pt]
%{\it If $x^*$  is a global minimizer of
%(TRS), there is a $\lambda^*\ge 0$ such that
%\begin{equation}
%(H+\lambda^*I)x^*=-c\label{trs}
%\end{equation}
%and $H+\lambda^*I$ is positive semidefinite. Moreover, assume
%$x^*$ and $\lambda^*\ge 0$ satisfy (\ref{trs}) and $H+\lambda^*I$ is positive semidefinite. If either $\|x^*\|=\Delta$ or $\|x^*\|<\Delta$ and $\lambda^*=0$, then $x^*$ is a global minimizer of (TRS).
%}
%\\[+5pt]
%
%In 1994, the local-nonglobal minimizer of (TRS)  was successfully
%characterized by Mart\'{\i}nez \cite{Mar}:
%  \\[+5pt]
%{\it Let $\alpha_1,~\alpha_2$ be the smallest two eigenvalues of $H$ and suppose $\alpha_1<\min\{0,\alpha_2\}$.
%Define the secular function
%
%and denote by $\phi'(\lambda)$ the derivative of  $\phi(\lambda)$. Then we have
%\begin{itemize}
%\item[(a)]
%If $x^*$  is a local-nonglobal minimizer of
%(TRS), then (\ref{trs}) holds with $\lambda^*\in (-\alpha_2,-\alpha_1)$, $\lambda^*\ge 0$ and $\phi'(\lambda^*)\ge 0$.
% \item[(b)] There exists at most one local-nonglobal minimizer of (TRS).
%  \item[(c)] If $\|x\|^*=\Delta$ holds for some $\lambda^*\in (\max\{-\alpha_2,0\},-\alpha_1)$ and $\phi'(\lambda^*)>0$,
%  then $x^*$ is a strict local minimizer of (TRS).
%\end{itemize}
%}
% \\[+5pt]

\section{Characterization of the Global Minimizers}

We first observe that the objective function $g(x)$ of (p-RS) is
coersive, i.e.,
\[
\lim_{\|x\|\rightarrow +\infty}g(x)=+\infty.
\]
Consequently, the global minimizer of (p-RS) always exists. The
starting point of the analysis is the first order and the second
order necessary conditions for any local minimizer of $g$.

%Due to the necessary optimality condition of a local minimizer, we have:
\begin{lem} \label{lem:1}
Assume that $\underline{x}$ is a local minimizer of (p-RS),~$p>2$.
It holds that
\begin{eqnarray}
&&\nabla g(\underline{x})=\left(H+\sigma\Vert
\underline{x}\Vert^{p-2} I\right)\underline{x}+c=0\label{con:1},\\
&&\nabla^2 g(\underline{x})=(H+\sigma\Vert \underline{x}\Vert^{p-2} I)+\sigma(p-2)\Vert
\underline{x}\Vert^{p-4}\underline{x}\underline{x}^{T}\succeq
0,\label{con:2}
\end{eqnarray}
where $\nabla g$, $\nabla^2 g$ denote the gradient and the Hessian
of $g(x)$, respectively.
\end{lem}

The next theorem shows that, a local minimizer $\underline{x}$
becomes global if and only if $H+\sigma\Vert
\underline{x}\Vert^{p-2} I\succeq0$. The necessity has been shown by
Theorem 2 in \cite{G10}. We only proves the sufficiency here.

\begin{thm} \label{thm:1}
The point $x^*$ is a global minimizer of (p-RS) for $p>2$ if and
only if it is a critical point satisfying $\nabla g(x^*)=0$ and
$H+\sigma \|x^*\|^{p-2}I\succeq0$. Moreover, the  $\ell_2$ norms of
all the global minimizers are equal.
\end{thm}
\proof If $x^*=0_n$, then $\sigma \|x^*\|^{p-2}=0$ so that
$c=-(H+\sigma \|x^*\|^{p-2}I)x^*=0$ and $H=H+\sigma
\|x^*\|^{p-2}I\succeq 0$. Consequently, $x^THx\ge0,~\forall x\in\Bbb
R^n.$ It follows that $x^*=0_n$ is a global minimizer since
\[
g(x)= \frac{1}{2}x^THx+c^Tx+\frac{\sigma}{p} \|x\|^p\ge \frac{\sigma}{p} \|x\|^p>0=g(0),~\forall x\neq 0_n=x^*.
\]
Now we assume $x^*\neq 0_n$, i.e., $\|x^*\|>0$. Define $Q=
H+\sigma\Vert x^*\Vert^{p-2}I$. According to the assumption,
$Q\succeq 0$. Then, for any $x\in \Bbb R^n$ and $x\neq x^*$, it
holds that
%\begin{eqnarray}
% g(x)&=&  \frac{1}{2}x^THx+c^Tx+\frac{\sigma}{p}\|x\|^p\nonumber \\
% &=& \frac{1}{2}x^TQx+c^Tx-
% \frac{1}{2}(\sigma\Vert x^*\Vert^{p-2})x^Tx+\frac{\sigma}{p}  \|x\|^p\nonumber\\
%  &=&\frac{1}{2}x^TQx+c^Tx+\frac{\sigma}{p}  \|x^*\|^p\left(
%   \left(\frac{\|x\|^2}{\|x^*\|^2}\right)^{\frac{p}{2}}-
% \frac{p}{2} \frac{\|x\|^2}{\|x^*\|^2}  \right)\nonumber\\
% &>&\frac{\sigma}{p}  \|x^*\|^p(1-\frac{p}{2})+ \frac{1}{2}x^TQx+c^Tx \label{fun:1}\\
%&\geq& \frac{\sigma}{p}  \|x^*\|^p -\frac{\sigma}{2}\|x^*\|^p +\frac{1}{2}x^{*T}Qx^*+c^Tx^*\label{fun:2}\\
%&=& g(x^*),\nonumber%\label{fun:3}
%\end{eqnarray}
\begin{eqnarray}
 g(x)&=&  \frac{1}{2}x^THx+c^Tx+\frac{\sigma}{p}\|x\|^p\nonumber \\
 &=& \frac{1}{2}x^TQx+c^Tx-
 \frac{1}{2}(\sigma\Vert x^*\Vert^{p-2})x^Tx+\frac{\sigma}{p}  \|x\|^p\nonumber\\
  &=&\frac{1}{2}x^TQx+c^Tx+\frac{\sigma}{p}  \|x^*\|^p\left(
   \left(\frac{\|x\|^2}{\|x^*\|^2}\right)^{\frac{p}{2}}-
 \frac{p}{2} \frac{\|x\|^2}{\|x^*\|^2}  \right)\label{fun:2}
\end{eqnarray}
Define $f(t)=t^{\frac{p}{2}},~p>2$. It is strictly convex for $t>0$.
Therefore,  \[ f(t)=t^{\frac{p}{2}} \ge
f(1)+f'(1)(t-1)=1+\frac{p}{2}(t-1),~\forall t>0.
\]
By substituting $t$ with $\frac{\|x\|^2}{\|x^*\|^2}$, we have
\begin{equation}\left(\frac{\|x\|^2}{\|x^*\|^2}\right)^{\frac{p}{2}}-
\frac{p}{2} \frac{\|x\|^2}{\|x^*\|^2}\ge 1-\frac{p}{2}.\nonumber
\end{equation}
Then,
\begin{eqnarray}
 g(x)&\ge&\frac{1}{2}x^TQx+c^Tx+\frac{\sigma}{p}  \|x^*\|^p(1-\frac{p}{2}).\label{fun:1}
\end{eqnarray}
By $Q\succeq0$, the lower bounding function of $g$ in the right hand
side of (\ref{fun:1}) is convex quadratic in terms of $x$. Since
$x^*$ satisfies $(H+\sigma\Vert x^*\Vert^{p-2}I)x^*=Qx^*=-c$, $x^*$
is a global minimizer of the convex function in the right hand side
of (\ref{fun:1}). As a consequence,
$$g(x)\ge\frac{1}{2}(x^*)^TQx^*+c^Tx^*+\frac{\sigma}{p}  \|x^*\|^p(1-\frac{p}{2})=g(x^*)$$
and $x^*$ is a global minimizer of (p-RS).

Finally, if $\|\hat x\|=\|x^*\|$, from (\ref{fun:2}) it can be seen
that
$$g(\hat x)=\frac{1}{2}(\hat x)^TQ{\hat x}+c^T{\hat x}+\frac{\sigma}{p}
\|x^*\|^p(1-\frac{p}{2}).$$ Then, $\hat x$ is also a global
minimizer of (p-RS) if and only if $\nabla g(\hat{x})=0$ and $\|\hat
x\|=\|x^*\|$.

 \hfill\eproof

To study the hidden convexity of (RS),
without loss of generality, we assume $H$ is diagonal, i.e.,
\begin{equation}
H={\rm Diag}(\alpha_1,\ldots,\alpha_n),~\alpha_1\le\ldots\le\alpha_n.\label{diag}
\end{equation}
Otherwise, let $H=U\Sigma U^T$ be the eigenvalue decomposition of $H$. Let $y=U^Tx$. Notice that $\|y\|=\|U^Tx\|=\|x\|$. We obtain a diagonal (RS) in terms of $y$.

\begin{pro}\label{lem:4}
Suppose $H$ is diagonal.
Let $x^{*}$ be the global minimizer of (RS), then we have
\begin{equation*}
c_ix_i^{*}\leq 0,~i=1,\ldots,n.
\end{equation*}
\end{pro}
\proof Let $\widetilde{x}=(-x_1^{*},
x_2^{*},x_3^{*},\ldots,x^{*}_n)$.  According to the definition of $x^*$, we have
$$
0\ge g(x^{*})-g(\widetilde{x})= c_1(x_1^{*}-\widetilde{x}_1)= 2c_1x_1^{*}.
$$
A similar argument applies for the other components. \hfill\eproof

Now we establish the hidden convexity of (RS).
According to  Proposition \ref{lem:4}, (RS) is equivalent to
\begin{eqnarray}\label{SDCDWPconstraint}
%\left\{
\begin{array}{ll}
\min\limits &  \sum\limits_{i=1}^n
\{\frac{\alpha_i}{2}x_i^2
+c_ix_i\}+\frac{\sigma}{p}\left( \sum\limits_{i=1}^n
 x_i^2\right)^{\frac{p}{2}}\\
\mbox{s.t.}&c_ix_i\leq0,\ \ i=1,\ldots,n.
\end{array}%\right.
\end{eqnarray}
Introducing the nonlinear one-to-one map:
\begin{equation}\label{x-y}
x_i=\left\{
\begin{array}{cl}\sqrt{z_i},&{\rm if}~c_i\leq 0,\\
-\sqrt{z_i},&{\rm if}~c_i> 0,
\end{array}\right.\  i = 1,\ldots,n,
\end{equation}
(RS) is equivalent to the following convex program:
\begin{eqnarray}\label{dualofdual}
\begin{array}{rll}
& \min\limits &
 -\sum\limits_{i=1}^n|c_i|\sqrt{z_i}
 +\frac{1}{2}\sum\limits_{i=1}^n\alpha_iz_i
+ \frac{\sigma}{p}\left(\sum\limits_{i=1}^nz_i\right)^{\frac{p}{2}}\\
&\mbox{s.t.}&z_i\geq 0,\ i = 1,\ldots,n.
\end{array}
\end{eqnarray}
Since (\ref{dualofdual}) is strictly convex when $p> 2$, again, we can see that $\sum\limits_{i=1}^nz^*_i=\|x^*\|^2$ is unique where $z^*$, $x^*$ are
any global minimizers of (\ref{dualofdual}) and (RS), respectively.

Before ending this section, we show that the set of the global minimizers of (RS), denoted by $\{x^*\}$,
is either a singleton or a $k$-dimensional sphere where $k$
is the multiplicity of the smallest
eigenvalue $\alpha_1$, i.e.,
\[
\alpha_1=\ldots=\alpha_k<\alpha_{k+1}\le\ldots\le\alpha_n.
\]
According to Theorem \ref{thm:1}, we first have
\[
\sigma\Vert x^*\Vert^{p-2}+\alpha_1\ge 0.
\]
\begin{itemize}
\item[$\bullet$] Suppose $c_1^2+\ldots+c_k^2>0$.
It follows from (\ref{con:1}) that
\[
\sigma\Vert x^*\Vert^{p-2}+\alpha_1> 0.
\]
Solving (\ref{con:1}) yields
\[
x^*_i=\frac{-c_i}{\sigma\Vert x^*\Vert^{p-2}+\alpha_i}, ~i=1,\ldots,n.
\]
By summing all $(x^*_i)^2$, we can derive that
\[
t^*=\Vert x^*\Vert^{p-2}
\]
is a nonnegative root of the following secular
function on a specific open interval:
\begin{equation}
 h(t)=\sum_{i=1}^n\frac{c_i^2}{(\sigma t+ \alpha_i)^2}-t^{\frac{2}{p-2}},
~t\in\left(- \frac{\alpha_1}{\sigma},+\infty\right). \label{h}
\end{equation}
Since
$\lim\limits_{t\rightarrow \max\{- \frac{\alpha_1}{\sigma},0\}}h(t)>0$,
$\lim\limits_{t\rightarrow+\infty}h(t)=-\infty$ and
$h(t)$ is strictly decreasing on $(- \frac{\alpha_1}{\sigma},+\infty)$,
the secular function $h(t)$ has a unique solution
$t^*$ on $(\max\{- \frac{\alpha_1}{\sigma},0\},+\infty)$. In this case, $x^*$ defined by
\begin{equation}
x_i^*=\frac{-c_i}{\sigma t^*+\alpha_i},~i=1,\ldots,n\label{x}
\end{equation}
is the unique global minimum solution of (RS).
\item[$\bullet$] Suppose $c_1^2+\ldots+c_k^2=0$. The secular
function (\ref{h}) reduces to
\begin{equation}
 h(t)=\sum_{i=k+1}^n\frac{c_i^2}{(\sigma t+ \alpha_i)^2}-t^{\frac{2}{p-2}},
~t\in\left(- \frac{\alpha_1}{\sigma},+\infty\right).  \label{h2}
\end{equation}
 There are two cases.
 \begin{itemize}
 \item[(1)] $\alpha_1>0$. Since $h(0)\ge 0$,  (\ref{h2}) has a unique nonnegative solution $t^*$. Then, $x^*$ satisfying (\ref{x})
is the unique global minimizer.
\item[(2)] $\alpha_1\le0$ and $h\left(- \frac{\alpha_1}{\sigma}\right)> 0$. Therefore, (\ref{h2}) has a unique nonnegative solution $t^*$ and hence $x^*$ satisfying (\ref{x})
is the unique global minimizer.
\item[(3)] $\alpha_1\le0$ and $h\left(- \frac{\alpha_1}{\sigma}\right)\le 0$. In this case, (\ref{h2}) has no solution. By Theorem \ref{thm:1}, any $x^*$ satisfying
    \begin{eqnarray}
    &&(x_1^*)^2+\ldots+(x_k^*)^2=-h\left(- \frac{\alpha_1}{\sigma}\right),
    \label{y1}\\
    &&x^*_i=-\frac{c_i}{\alpha_i-\alpha_1},  ~i=k+1,\ldots,n \label{y2}
    \end{eqnarray}
    is a global minimizer. Namely, the global minimum
    solution set forms a $k$-dimensional sphere
    centered at
    $(0,\cdots,0,-\frac{c_{k+1}}{\alpha_{k+1}-\alpha_1},\cdots,-\frac{c_{n}}{\sigma_{n}-\sigma_1})$
    with the radius
    $\sqrt{\left(\frac{\alpha_1}{\sigma}\right)^{\frac{2}{p-2}}-\sum_{i=k+1}^n\frac{c_i^2}{(\alpha_i-\alpha_1)^2}}$.
    \end{itemize}
\end{itemize}

\section{Characterization of the Local-Nonglobal Minimizer}
In this section, we establish the necessary and sufficient optimality condition for the local-nonglobal minimizer of (RS).

Let $\alpha_1\le\ldots\le\alpha_n$ be the eigenvalues of $H$. Throughout this section, %without loss of generality, we always assume $H$ is diagonal as in (\ref{diag}).
we assume $\alpha_1<0$. That is,
$H\not\succeq 0$. Otherwise, (RS) is a convex minimization problem and hence has no local-nonglobal minimizer.

\begin{lem}
Suppose $H\not\succeq 0.$ Then $0_n$ is not a local minimizer of (RS).
\end{lem}
\proof
Suppose $0_n$ is a local minimizer of (RS). Then the necessary optimality conditions
(\ref{con:1})-(\ref{con:2}) imply that
\[
c=0,~H\succeq 0,
\]
which is contradiction.
\endproof

\begin{lem} \label{lem:2}
Suppose $n\geq 2$. Let $\underline{x}$ be a local minimizer of (RS). It holds that
\begin{eqnarray}
\sigma\Vert \underline{x}\Vert^{p-2}+\alpha_2\geq 0.
\label{lem:eq}
\end{eqnarray}
Furthermore, if $\alpha_1<\alpha_2$, then
\begin{eqnarray} \label{lem:ineq}
\sigma\Vert \underline{x}\Vert^{p-2}+\alpha_2> 0.
\end{eqnarray}
\end{lem}
\proof Without loss of generality, we can assume $H$ is a diagonal matrix, i.e.,
$H={\rm Diag}(\alpha_1,\ldots,\alpha_n)$. Suppose the statement (\ref{lem:eq}) is not true, then
$\sigma\Vert \underline{x}\Vert^{p-2}+\alpha_2<0$.
\[
\sigma\Vert \underline{x}\Vert^{p-2}+\alpha_1\le \sigma\Vert \underline{x}\Vert^{p-2}+\alpha_2<0.
\]
 Let $e_1$ and $e_2$ be the first two columns of $I$, respectively. We consider the following two cases.
\begin{itemize}
\item[(a)] Suppose $\underline{x}_1=e_1^T\underline{x}=0$. It follows from the necessary condition (\ref{con:2}) that
\[
0\leq e_1^T(\sigma(p-2)\Vert
\underline{x}\Vert^{p-4}\underline{x}\underline{x}^{T}+ \sigma\Vert \underline{x}\Vert^{p-2} I+H)e_1
=\sigma\Vert \underline{x}\Vert^{p-2}+\alpha_1<0,
\]
which is a contradiction.
\item[(b)]Suppose $\underline{x}_1=e_1^T\underline{x}\neq0$.  It follows from the necessary condition  (\ref{con:2}) that
\begin{eqnarray*}
0&\leq
&((-\underline{x}_2)e_1+(\underline{x}_1)e_2)^T(\sigma(p-2)\Vert
\underline{x}\Vert^{p-4}\underline{x}\underline{x}^{T}\\
&~~&+ \sigma\Vert \underline{x}\Vert^{p-2} I+H)((-\underline{x}_2)e_1+(\underline{x}_1)e_2)\\
&=&(\sigma\Vert \underline{x}\Vert^{p-2}+\alpha_1)(\underline{x}_2)^2+
(\sigma\Vert \underline{x}\Vert^{p-2}+\alpha_2)(\underline{x}_1)^2<0,
\end{eqnarray*}
which is a contradiction.
\end{itemize}
Therefore, the statement (\ref{lem:eq}) holds true.

Now we assume $\alpha_1<\alpha_2$ and suppose that the statement
(\ref{lem:ineq}) is not true. Then we have
\begin{eqnarray}\label{lem:eq_1}
\sigma\Vert \underline{x}\Vert^{p-2}+\alpha_2=0,
\end{eqnarray}
%By (\ref{lem:eq_1}), we know that the second order necessary
with which the necessary optimality
condition (\ref{con:2}) becomes
\begin{eqnarray}
0&\preceq& \sigma(p-2)\Vert
\underline{x}\Vert^{p-4}\underline{x}\underline{x}^{T}+ \sigma\Vert \underline{x}\Vert^{p-2} I+H \nonumber\\
&=& \sigma(p-2)\Vert
\underline{x}\Vert^{p-4}\underline{x}\underline{x}^{T}+H-
\alpha_2I.\label{x_2}
\end{eqnarray}
Consequently, the first two leading principal minors of the matrix in
(\ref{x_2}) are nonnegative, i.e.,
\begin{eqnarray}\label{x_0}
\sigma(p-2)\Vert
\underline{x}\Vert^{p-4} \underline{x}_1^{2}+\alpha_1-\alpha_2\ge0
\end{eqnarray}
and
\begin{eqnarray}
&&\det\left\{\sigma(p-2)\Vert
\underline{x}\Vert^{p-4}\left[\begin{array}{cc} \underline{x}_1^{2}&
\underline{x}_1\underline{x}_2
\\ \underline{x}_1\underline{x}_2& \underline{x}_2^{2}
\end{array}\right]+\left[\begin{array}{cc}\alpha_1-\alpha_2&0\\0&0 \end{array}\right]
\right\} \nonumber\\
&&= \sigma(p-2)\Vert
\underline{x}\Vert^{p-4}(\alpha_1-\alpha_2)\underline{x}_2^{2}\ge 0.\label{x1}
\end{eqnarray}
Since $\alpha_1-\alpha_2<0$, the inequalities (\ref{x_0}) and (\ref{x1}) imply
$\underline{x}_1\not=0$ and $\underline{x}_2=0$, respectively. Then it follows from the necessary optimality condition (\ref{con:1}) that
obtain that $c_2=0$ and
$$\underline{x}_1=\frac{-c_1}{\sigma\Vert
\underline{x}\Vert^{p-2}+\alpha_1}
=\frac{ c_1}{\alpha_2-\alpha_1}.
$$
Without loss of generality, we assume that $c_1>0$, which implies that
$\underline{x}_1>0$. Then, according to (\ref{lem:eq_1}) and the
fact $\underline{x}_2=0$, we have
$$\underline{x}_1=\sqrt{\left(\frac{-\alpha_2}{\sigma}\right)^{\frac{2}{p-2}}-\sum_{i=3}^n\underline{x}_i^{2}}.$$
Consider the following parametric curve in $\mathbb{R}^n$:
\begin{eqnarray}\label{curve}
\gamma(t)&=& \{(k(t),t,\underline{x}_3,\ldots,\underline{x}_n)\vert~ \nonumber\\
&& k(t)=\sqrt{
\left(\frac{-\alpha_2}{\sigma}\right)^{\frac{2}{p-2}}-t^2-\sum_{i=3}^n
\underline{x}_i^{2}}=\sqrt{\underline{x}_1^{2}-t^2},
~t\in \Bbb R \}
\end{eqnarray}
where $\gamma(0)=\gamma(\underline{x}_2)=\underline{x}$, i.e.,
$\gamma(t)$ passes through $\underline{x}$ at $t=0$. Evaluating
$g(x)$ on $\gamma(t)$, we have
\begin{eqnarray*}
&&g(\gamma(t))\\
&=&\frac{\sigma}{p}\left( k(t)^2 + t^2
+ \sum_{i=3}^n\underline{x}_i^{2}
\right)^{\frac{p}{2}} + \frac{\alpha_1}{2}k(t)^2+\frac{\alpha_2}{2}t^2+\sum\limits_{i=3}^n\frac{\alpha_i}{2}
\underline{x}_i^{2}
+c_1k(t)+\sum_{i=3}^nc_i\underline{x}_i\\
&=&\frac{\sigma}{p}\left(\underline{x}_1^{2} + \sum_{i=3}^n\underline{x}_i^{2}
\right)^{\frac{p}{2}}+
 \frac{\alpha_1}{2}\underline{x}_1^{2}+
\sum_{i=3}^n\frac{\alpha_i}{2}\underline{x}_i^{2}
+\frac{\alpha_2-\alpha_1}{2}t^2+c_1\sqrt{\underline{x}_1^{2}-
t^2}+\sum_{i=3}^nc_i\underline{x}_i.
\end{eqnarray*}
Since $\underline{x}$ is a local minimizer of $g(x)$,
$t=0$ must be a local minimum point of $g(\gamma(t))$.
However, this conclusion contradicts to the fact that
$$\frac{d}{dt}g(\gamma(0))=\frac{d^2}{dt^2}g(\gamma(0))=\frac{d^3}{dt^3}g(\gamma(0))=0,
~\frac{d^4}{dt^4}g(\gamma(0))=-\frac{3(\alpha_2-\alpha_1)}{\underline{x}_1^{2}}<0.$$
Consequently, the statement (\ref{lem:ineq}) holds true under the additional assumption
$\alpha_1<\alpha_2$.
\hfill\eproof

As the main result in this section, we
 establish the necessary and sufficient condition
for local-nonglobal minimizer of (RS).
\begin{thm}\label{thm:2}
$\underline{x}$ is a local-nonglobal minimizer of (RS) if and only if
\begin{equation}
\underline{x}
=-
\left( \sigma\underline{t}^* I+
H\right)^{-1}c,\label{tw}
\end{equation}
%\underline{x}
%=\left(\frac{-c_1}{\sigma\underline{t}^*+\alpha_1},\ldots,\frac{-c_n}
%{\sigma\underline{t}^* +\alpha_n} \right),
where $\underline{t}^*$ is a root of the secular function
\begin{equation}
 h(t)=\|\left( \sigma t  I+
H\right)^{-1}c\|^2-t^{\frac{2}{p-2}},
 ~~t\in\left(\max\left\{- \frac{\alpha_2}{\sigma},0\right\}, - \frac{\alpha_1}{\sigma}\right)\label{varphi0}
\end{equation}
%\begin{equation}
% h(t)=\sum_{i=1}^n\frac{c_i^2}{(\sigma t+ \alpha_i)^2}-t^{\frac{2}{p-2}},
% ~~t\in\left(\max\{- \frac{\alpha_2}{\sigma},0\}, - \frac{\alpha_1}{\sigma}\right),\label{varphi}
%\end{equation}
such that $h'(\underline{t}^*)>0$.
\end{thm}
\proof Without loss of generality, we can assume $H$ is a diagonal matrix, i.e.,
$H={\rm Diag}(\alpha_1,\ldots,\alpha_n)$.
It is sufficient to consider the nontrivial case  $n\geq 2$,
since for $n=1$,  we will see that it amounts to setting $\alpha_2=\infty$ in the following proof.

According to Lemma \ref{lem:2} and Theorem \ref{thm:1}, the local-nonglobal minimizer
$\underline{x}$ of  (RS) exits only if
\begin{equation}
-\alpha_2< \sigma\Vert \underline{x}\Vert^{p-2}<-\alpha_1.\label{loc:1}
\end{equation}
It follows that the diagonal
matrix $\sigma\Vert \underline{x}\Vert^{p-2} I+H$ is
nonsingular with its first diagonal element being negative and
others positive. Solving (\ref{con:1}), we obtain
\begin{equation}
\underline{x}_i=\frac{-c_i}{\sigma\Vert \underline{x}\Vert^{p-2}+\alpha_i}, ~~i=1,\ldots,n.\label{sol}
\end{equation}
The necessary optimality condition (\ref{con:2}) implies that
\[
\sigma(p-2)\Vert
\underline{x}\Vert^{p-4}\underline{x}^{2}_1+ \sigma\Vert \underline{x}\Vert^{p-2}+\alpha_1\geq 0.
\]
Then it follows from the right hand side of (\ref{loc:1}) that
\begin{equation}
\|x\|>0,\label{loc:2}
\end{equation}
%\begin{equation}\label{a1}
and moreover,
\begin{equation}
\underline{x}_1\not=0,\label{x10}
\end{equation}
%\end{equation}
%and hence $c_1\not=0$ according to (\ref{sol}).
Putting all $\underline{x}_i$ in (\ref{sol}) together yields
\begin{equation}
\sum_{i=1}^n\frac{c_i^2}{(\sigma\Vert \underline{x}\Vert^{p-2}+\alpha_i)^2}=\Vert \underline{x}\Vert^2.\label{loc:3}
\end{equation}
As a summary of (\ref{loc:1}), (\ref{loc:2}) and (\ref{loc:3}),
\[
t^*=\Vert
\underline{x}\Vert^{p-2}
\]
is a root of the following secular
function on a specific open interval:
\begin{equation}
 h(t)=\sum_{i=1}^n\frac{c_i^2}{(\sigma t+ \alpha_i)^2}-t^{\frac{2}{p-2}},
 ~~t\in\left(\max\left\{- \frac{\alpha_2}{\sigma},0\right\}, - \frac{\alpha_1}{\sigma}\right),\label{varphi}
\end{equation}
which is the diagonal version of (\ref{varphi}).
Notice that each root of $h(t)=0$ can only correspond to one local-nonglobal minimizer of (RS) due to (\ref{sol}).
Taking a simple calculation of (\ref{varphi}), we have
\begin{eqnarray}
%\varphi(t)&=&\sum_{i=1}^n\frac{b_i^2}{(t-2\alpha+\lambda_i)^2}-t,\label{varphi}\\
h'(t) = -\sum_{i=1}^n\frac{2\sigma c_i^2}{(\sigma t+\alpha_i)^3}-\frac{2}{p-2}t^{\frac{4-p}{p-2}}.\label{a4}
\end{eqnarray}
%Since $\sigma\Vert \underline{x}\Vert^{p-2}I+H$ is
%nonsingular,
We notice that the necessary optimality condition (\ref{con:2}) is equivalent to
\begin{equation}
\sigma(p-2)\Vert
\underline{x}\Vert^{p-4}(\Gamma\underline{x})(\Gamma\underline{x})^T+ {\rm
Diag}(-1,1,\ldots,1)  \succeq 0,\label{con:t}
\end{equation}
where
\begin{equation}
\Gamma={\rm Diag}\left(\frac{1}{\sqrt{-\sigma\Vert \underline{x}\Vert^{p-2}-\alpha_1}},  \frac{1}{\sqrt{\sigma\Vert \underline{x}\Vert^{p-2}+ \alpha_2}},\ldots, \frac{1}{\sqrt{\sigma\Vert \underline{x}\Vert^{p-2}+\alpha_n}}\right). \label{gamma:1}
\end{equation}
Since the determinant of the positive semidefinite matrix in (\ref{con:t}) is nonnegative, we have
\begin{eqnarray}
0&\leq& {\rm det}\left(\sigma(p-2)\Vert
\underline{x}\Vert^{p-4}(\Gamma \underline{x}) (\Gamma
\underline{x})^T+
{\rm Diag}(-1,1,\ldots,1)\right)\nonumber \\
&=&{\rm det}({\rm Diag}(-1,1,\ldots,1))\times \nonumber\\
&~~&{\rm det}\left(\sigma(p-2)\Vert
\underline{x}\Vert^{p-4}{\rm
Diag}(-1,1,\ldots,1)
(\Gamma \underline{x}) (\Gamma \underline{x})^T+I\right)\nonumber \\
&=&-1\times\left(\sigma(p-2)\Vert
\underline{x}\Vert^{p-4}(\Gamma\underline{x})^T{\rm Diag}(-1,1,\ldots,1)(\Gamma\underline{x})+1\right)\nonumber \\
&=&-\sum_{i=1}^n\frac{\sigma(p-2)\Vert
\underline{x}\Vert^{p-4}c_i^2}{(\sigma\Vert \underline{x}\Vert^{p-2}+\alpha_i)^3}-1\nonumber\\
&=& (\frac{p}{2}-1)\Vert
\underline{x}\Vert^{p-4}h'(\Vert \underline{x}\Vert^{p-2})\nonumber\\
&=& (\frac{p}{2}-1)\Vert
\underline{x}\Vert^{p-4}h'(t^*). \nonumber
\end{eqnarray}
It follows from $p>2$ and (\ref{loc:2}) that  $h'(t^*)\ge 0$. Now,
it remains to show that $ h'(\underline{t}^*)>0$. Suppose this is not true, we have $h'(\underline{t}^*)=0$. Therefore,
we obtain
\begin{eqnarray}
&&\det\left(\sigma(p-2)\Vert
\underline{x}\Vert^{p-4}\underline{x}\underline{x}^{T}+\sigma\Vert
\underline{x}\Vert^{p-2}I+H\right)\nonumber\\
&=&\frac{{\rm det}\left(\sigma(p-2)\Vert
\underline{x}\Vert^{p-4}(\Gamma \underline{x}) (\Gamma
\underline{x})^T+
{\rm Diag}(-1,1,\ldots,1)\right)}{{\det}^2(\Gamma)}\nonumber\\
&=&\frac{(\frac{p}{2}-1)\Vert
\underline{x}\Vert^{p-4}h'(\underline{t}^*)}{ {\det}^2(\Gamma)}\label{determinant}\\
&=&0 \nonumber
\end{eqnarray}
and thus there is a $u=(u_1,\ldots,u_n)^T\neq 0$ such that
\begin{equation}
\sigma(p-2)\Vert
\underline{x}\Vert^{p-4}\underline{x}\underline{x}^{T}u+\left(\sigma\Vert
\underline{x}\Vert^{p-2}I+H \right)u=0,\label{new:000}
\end{equation}
%From (\ref{new:000}), we can write
or equivalently,
\begin{equation*}
u_i=\frac{-\sigma(p-2)\Vert
\underline{x}\Vert^{p-4}\underline{x}_i(u^T\underline{x})}{\sigma\Vert
\underline{x}\Vert^{p-2} +\alpha_i},~i=1,2,\ldots,n.
\end{equation*}
Since $u\neq 0$,  it holds  that
\begin{equation}
u^T\underline{x}\neq 0.\label{nzero}
\end{equation}
Define
\[
q(\beta):=g(\underline{x}+\beta u).
\]
We can verify that
\begin{eqnarray*}
q'(\beta)&=& \nabla g(\underline{x}+\beta u)u, \\
q''(\beta)&=& u^T \nabla^2g(\underline{x}+\beta u)u,\\
%q'''(\beta)&=&2\sigma(p-2)\Vert
%\underline{x}+\beta u\Vert^{p-4}  u^T(\underline{x}+\beta u)u^Tu\\
%&~~&+\sigma(p-2)(p-4)\Vert
%\underline{x}+\beta u\Vert^{p-6}u^Tu\beta
%+ \sigma(p-2)\Vert
%\underline{x}+\beta u\Vert^{p-4}(u^Tu)^2\beta.
q'''(\beta)&=&3\sigma(p-2)\Vert
\underline{x}+\beta u\Vert^{p-4} (u^T\underline{x}+\beta u^Tu)u^Tu\\
&~~&+\sigma(p-2)(p-4)\Vert
\underline{x}+\beta u\Vert^{p-6}(u^T\underline{x}+\beta u^Tu)^3.
\end{eqnarray*}
The necessary optimality condition (\ref{con:1}) implies that
$q'(0)=0$. According to the definition of $u$, we have
$q''(0)=0$. However, (\ref{nzero}) implies that
\begin{eqnarray*}
\left(q'''(0)\right)^2&=&\sigma^2(p-2)^2\Vert
\underline{x}\Vert^{2(p-6)}(u^T\underline{x})^6\left( 3\frac{(\underline{x}^T\underline{x})  (u^Tu)}{(u^T\underline{x})^2}
+ (p-4)\right)^2\\
&\ge&\sigma^2(p-2)^2\Vert
\underline{x}\Vert^{2(p-6)}(u^T\underline{x})^6 (p-1)^2\\
&>&0,
\end{eqnarray*}
where the first inequality follows from Cauchy-Schwartz inequality.
It contradicts to the fact that $\underline{x}$ is a local
minimizer of (RS). Therefore, $
h'(\underline{t}^*)>0$ and the necessary proof is complete.

It remains for us to give the sufficient proof.
Let $t^*\in(\max\{- \frac{\alpha_2}{\sigma},0\}, - \frac{\alpha_1}{\sigma})$ be a root of the secular function (\ref{varphi}) such that $h'(t^*)>0$. Define
$\underline{x}$ as in (\ref{tw}). Then we have
\[
\Vert\underline{x}\Vert^{2}=\sum_{i=1}^n\frac{c_i^2}{(\sigma\underline{t}^*+\alpha_i)^2}
=(\underline{t}^*)^{\frac{2}{p-2}},
\]
%Suppose there is a $\underline{t}^*\in(\max\{- \frac{\alpha_2}{\sigma},0\}, - \frac{\alpha_1}{\sigma})$ such that
%\[
%h(\underline{t}^*)=\sum_{i=1}^n\frac{c_i^2}{(\sigma\underline{t}^*+\alpha_i)^2}
%-(\underline{t}^*)^{\frac{2}{p-2}}=0
%\]
%and $h'(\underline{t}^*)>0$.
that is,
$\underline{t}^*=\Vert\underline{x}\Vert^{p-2}$.
Consequently, $\underline{x}$
satisfies the first-order necessary optimality condition (\ref{con:1}).
Moreover, the diagonal matrix $\sigma\Vert
\underline{x}\Vert^{p-2}I+H$ is nonsingular with positive diagonal
elements  except for the first one. By Weyl's inequality (see
\cite{Horn}, Theorem 4.3.1), we have
\begin{eqnarray}
\lambda_i\left(\nabla^2 g(\underline{x})\right)&=&\lambda_i\left(\sigma(p-2)\Vert
\underline{x}\Vert^{p-4}\underline{x}\underline{x}^{T}+\sigma\Vert
\underline{x}\Vert^{p-2}I+H\right) \nonumber\\
&\ge&
\lambda_1\left(\sigma(p-2)\Vert
\underline{x}\Vert^{p-4}\underline{x}\underline{x}^{T}\right)+\lambda_i
\left( \sigma\Vert
\underline{x}\Vert^{p-2}I+H\right)\nonumber\\
&=&\lambda_i\left(\sigma\Vert
\underline{x}\Vert^{p-2}I+H\right)\nonumber\\
&>&0,~~~~~\hbox{for}~i=2,3,\ldots,n,\label{new:1}
\end{eqnarray}
where $\lambda_i(P)$ is the $i$th smallest eigenvalue of $P$.
Since $h'(\underline{t}^*)>0$, by (\ref{determinant}), we have
\begin{eqnarray}
\prod_{i=1}^n\lambda_i\left(\nabla^2 g(\underline{x})\right)&=&\det\left(\nabla^2 g(\underline{x})\right)\nonumber\\
&=&\det\left(\sigma(p-2)\Vert
\underline{x}\Vert^{p-4}\underline{x}\underline{x}^{T}+\sigma\Vert
\underline{x}\Vert^{p-2}I+H\right)\nonumber\\
&=&\frac{(\frac{p}{2}-1)\Vert
\underline{x}\Vert^{p-4}h'(\underline{t}^*)}{{\det}^2(\Gamma)}\nonumber\\
&>&0,\label{new:2}
\end{eqnarray}
%where $\Gamma$ is defined in (\ref{gamma:1}).
Combining (\ref{new:1}) with (\ref{new:2}), we have
\[
\lambda_1\left(\nabla^2 g(\underline{x})\right)>0,
\]
or equivalently,
\[
\nabla^2 g(\underline{x})=\sigma(p-2)\Vert
\underline{x}\Vert^{p-4}\underline{x}\underline{x}^{T}+\sigma\Vert
\underline{x}\Vert^{p-2}I+H\succ0.
\]
This is a sufficient condition
to guarantee that $\underline{x}$ is a local minimizer of  (RS).
The proof is complete.
 \hfill\eproof

Theorem \ref{thm:2} and its proof  provide some simple sufficient conditions for
having no local-nonglobal minimizer.
\begin{cor}\label{cor:1}
When one of the following conditions is met:
\begin{itemize}
\item[(a)]
$\alpha_1\geq 0$;
\item[(b)] $\alpha_1=\alpha_2$;
\item[(c)] $v^Tc=0$, where $v$ is the eigenvector of $H$ corresponding to $\alpha_1$;
\end{itemize}
any local minimizer of  (RS) is globally optimal.
\end{cor}
\proof
In Case (a), $g(x)$ is convex and hence any local minimizer is globally optimal. In Case (b), it is trivial to see that the secular function (\ref{varphi})  has no solution. Therefore, according to Theorem \ref{thm:2}, the local-nonglobal minimizer does not exist.
Suppose (RS) has a local-nonglobal minimizer in Case (c). Let $H=U\Sigma U^T$ be the eigenvalue decomposition of $H$. Introducing $y=U^Tx$, we obtain a diagonal version of (RS) with respect to $y$:
\[
\min_{y\in \Bbb R^n} \left\{
g(x)= \frac{1}{2}y^T\Sigma x+\widetilde{c}^Ty+\frac{\sigma}{p} \|y\|^p \right\}
\]
where $\widetilde{c}=U^Tc$.
According to (\ref{sol}) and (\ref{x10}) in the necessary proof of Theorem \ref{thm:2}, a necessary condition for the secular function (\ref{varphi}) having a solution is that $\widetilde{c}_1\neq 0$. We obtain a contradiction by noting that $\widetilde{c}_1=(U^Tc)_1=v^Tc$, where $v$ is the eigenvector of $H$ corresponding to $\alpha_1$.
  \hfill\eproof

The second corollary of Theorem \ref{thm:2} can be regarded as the similar version of Proposition \ref{lem:4} for the local-nonglobal minimizer.
\begin{cor}\label{cor:2}
Suppose $H$ is diagonal.
Let $\underline{x}$ be the local-nonglobal minimizer of (RS), then we have
\begin{equation}
c_1\underline{x}_1> 0,~c_i\underline{x}_i \leq 0,~i=2,3,\ldots,n.
\end{equation}
\end{cor}
\proof
Following (\ref{tw}) and (\ref{varphi0}), we immediately have
\[
c_1\underline{x}_1 \ge 0,~c_i\underline{x}_i \leq 0,~i=2,3,\ldots,n.
\]
The fact $\underline{x}_1\neq 0$ is shown in (\ref{x10}) and
the statement $c_1\neq 0$ follows from  (\ref{x10}) and (\ref{sol}).
  \hfill\eproof

As an application of Corollary \ref{cor:2}, similar to (\ref{dualofdual})
we see that finding the local-non-global minimizer of
(RS) is actually equivalent to globally minimizing the following nonconvex program:
\begin{eqnarray}
\begin{array}{rll}
& \min\limits &
 |c_i|\sqrt{z_i}-\sum\limits_{i=2}^n|c_i|\sqrt{z_i}
 +\frac{1}{2}\sum\limits_{i=1}^n\alpha_iz_i
+ \frac{\sigma}{p}\left(\sum\limits_{i=1}^nz_i\right)^{\frac{p}{2}}\\
&\mbox{s.t.}&z_i\geq 0,\ i = 1,\ldots,n.
\end{array}\nonumber
\end{eqnarray}

As the last corollary of Theorem \ref{thm:2}, we have
\begin{thm}\label{thm:3}
(RS) with $p>2$ has at most one local-nonglobal minimizer.
\end{thm}
\proof
First we observe that the secular function (\ref{varphi0}) %becomes (\ref{varphi}).
has the same roots as
\[
 p(t)=\log\left(\|\left( \sigma t  I+
H\right)^{-1}c\|^2\right)-\frac{2}{p-2}\log(t),
~t\in\left(\max\left\{- \frac{\alpha_2}{\sigma},0\right\}, - \frac{\alpha_1}{\sigma}\right).
\]
Without loss of generality, we assume $H$ is diagonal.  If $c_1=0$, then (RS) has no local-nonglobal minimizer according to Corollary \ref{cor:2}. So, we assume $c_1\neq 0$.
Then, we have
\[
p''(t)=\frac{\sum_{i=1}^n\frac{6\sigma^2 c_i^2}{(\sigma t+\alpha_i)^4}}{\sum_{i=1}^n\frac{c_i^2}{(\sigma t+\alpha_i)^2}}
-\frac{\left(\sum_{i=1}^n\frac{2\sigma c_i^2}{(\sigma t+\alpha_i)^3}\right)^2}{\left(\sum_{i=1}^n\frac{c_i^2}{(\sigma t+\alpha_i)^2}\right)^2}.
\]
Define two vectors in $\Bbb R^n$:
\[%\begin{eqnarray*}
a=\left(\frac{\sqrt{6} \sigma c_1}{(\sigma t+\alpha_1)^2},\ldots,
\frac{\sqrt{6} \sigma c_n}{(\sigma t+\alpha_n)^2}\right)^T,~
b=\left(\frac{ c_1}{\sigma t+\alpha_1 },\ldots,
\frac{ c_n}{ \sigma t+\alpha_n }\right)^T.
\]%\end{eqnarray*}
Applying Cauchy-Schwartz inequality, we obtain
\begin{eqnarray*}
\left(\sum_{i=1}^n\frac{2\sigma c_i^2}{(\sigma t+\alpha_i)^3}\right)^2
&<& (a^Tb)^2\\
&\le& (a^Ta)(b^Tb)\\
&=&\left(\sum_{i=1}^n\frac{6\sigma^2 c_i^2}{(\sigma t+\alpha_i)^4}\right)\left(\sum_{i=1}^n\frac{c_i^2}{(\sigma t+\alpha_i)^2}\right).
\end{eqnarray*}
Therefore, $p''(t)>0$ for all $t$ such that $p(t)$ is well-defined. It follows that $p(t)$ is strictly convex for $t\in\left(\max\left\{- \frac{\alpha_2}{\sigma},0\right\}, - \frac{\alpha_1}{\sigma}\right)$. Thus,
$p(t)$, as well as $h(t)$, has at most two real roots in this interval. Let $t_1<t_2$ be the only two roots of $h(t)$.
Suppose  $h'(t_1)>0$ and $h'(t_{2})>0$.
Then, for sufficiently small $\epsilon\in(0,\frac{t_{2}-t_1}{2})$, we have
\[
h(t_1+\epsilon)>h(t_1)=0,~h(t_{2}-\epsilon)<h(t_{2})=0.
\]
Therefore, there is a $\widetilde{t}\in[t_1+\epsilon,t_{2}-\epsilon]$ such that $h(\widetilde{t})=0$, which is a contradiction.
consequently, the secular function $h(t)$ has at most one real root
satisfying $h'(t) > 0$. Following Theorem \ref{thm:2},
the proof is complete.
\hfill\eproof

\section{(RS) with linear inequality constraints}
In this section, we study  $({\rm RS}_m)$ (\ref{RSm:1})-(\ref{RSm:2}).
For a special case $p=4$,
we first show $({\rm RS}_m)$ is  NP-hard when $m>n$. Then, as an application
of Theorems \ref{thm:2} and \ref{thm:3}, we show $({\rm RS}_m)$ can be solved
in polynomial time when $m$ is a fixed number.
\subsection{NP-hardness}
Let
$S=\{x\in[0,1]^n\mid~ Ax\leq b\}$ be nonempty, where
$A\in\Bbb R^{m\times n}$ and $b\in\Bbb R^m$ are with integer elements, and
for $i=1,\ldots,n$, the $i$-th row of
$A$, denoted by $a_i$, satisfies that  $\|a_i\|_1\geq 2$.

\begin{lem}[\cite{Hs}]\label{lemma1}
For any vertex $x=(x_1,\cdots, x_n)^T$ of the polytope $S$, if
$x\not\in\{0, 1\}^n$, then it holds that
\begin{equation}\label{eq:1}
x^T(e-x)\geq  \left\{\begin{array}{ll}\frac{\max_{1\leq j\leq m}\|a_j\|_{\infty}-1}{\max_{1\leq j\leq m}\|a_j\|_{\infty}^2},& {\rm if}~ \max_{1\leq j\leq m}\|a_j\|_{\infty}\geq 2,\\
\frac{1}{2},& {\rm if}~\max_{1\leq j\leq
m}\|a_j\|_{\infty}=1,\end{array}\right.
\end{equation}
where $e$ is a vector of dimension $n$ with all components equal to one.
%\begin{equation}
%g(A)= \left\{\begin{array}{ll}\frac{\max_{1\leq j\leq m}\|a_j\|_{\infty}^2}{\max_{1\leq j\leq m}\|a_j\|_{\infty}-1},& {\rm if}~ \max_{1\leq j\leq m}\|a_j\|_{\infty}\geq 2,\\
%2,& {\rm if}~\max_{1\leq j\leq
%m}\|a_j\|_{\infty}=1.\end{array}\right.\label{gA}
%\end{equation}
\end{lem}
Now, we consider the following $k$-dispersion-sum problem:
\begin{eqnarray}
{\rm(KDSP)}~~d^*=&\min &x^T(-D)x \label{p2:1}\\
&{\rm s.t.}&e^Tx=k,~x\in \{0,1\}^n.\label{p2:2}
\end{eqnarray}
It is to locate $k$ facilities at some of $n$
predefined locations by maximizing the distance sum between the $k$ established facilities, where the distance between two facilities $i$ and $j$ is given by
a square matrix $D = (d_{ij})$, $i, j = 1, 2, \ldots, n$.
(KDSP) is NP-hard, even if the distance matrix satisfies the triangle
inequality, see \cite{E90,Ha}.

Define the continuous relaxation of (KDSP) as
\[
d^c= \min_{e^Tx=k,~x\in [0,1]^n}x^T(-D)x.
\]
It trivially holds that
$
d^c\le d^*.
$
For any $\theta\ge 4(d^*-d^c)$, we obtain
\begin{eqnarray*}
&&\min_{e^Tx=k,~x\in [0,1]^n}\left\{d(x):=-x^TDx+\theta \left(e^Tx-x^Tx\right)^2 \right\}\\
&&= \min\left\{\min_{e^Tx=k,~x\in [0,1]^n\setminus \{0,1\}^n}d(x),~\min_{e^Tx=k,~x\in \{0,1\}^n}d(x)\right\}\\
&&\ge \min\left\{\min_{e^Tx=k,~x\in [0,1]^n}-x^TDx+\min_{e^Tx=k,~x\in [0,1]^n\setminus \{0,1\}^n}\theta\left(e^Tx-x^Tx\right)^2,~d^*\right\}\\
&&\ge \min\left\{ d^c+ \frac{1}{4}\theta,~d^*\right\}\\
&&=d^*,
\end{eqnarray*}
where the second inequality holds since
 \[
 x^T(e-x)\geq \frac{1}{2},~\forall x\in [0,1]^n\setminus \{0,1\}^n,~e^Tx=k,
 \]
which follows from Lemma \ref{lemma1}.
Therefore,
(KDSP) has been reduced to the following special case of $({\rm RS}_{n+1})$ with $p=4$:
\begin{eqnarray}
&\min &-x^TDx+\theta \left(k-x^Tx\right)^2=
-x^T(D+2\theta k\cdot I)x+\theta \|x\|^4 +\theta k^2   \\
&{\rm s.t.}&e^Tx=k,~x\in [0,1]^n.
\end{eqnarray}
As a summary, we have the following result:
\begin{thm}
When $p=4$, $({\rm RS}_{m})$ with $m> n$ is NP-hard.
\end{thm}
\subsection{Polynomially Solvable Cases}

Consider $({\rm RS}_{m})$ with $p=4$ and $m$ being a fixed number.
The approach applied in this subsection inherits from \cite{HS}.

Let $X_0^*$ denote the set of the global minimizers of $({\rm RS})$.
According to the discussion at the end of Section 2,  $X_0^*$ is either
a singleton or a $k$-dimensional sphere, both can be obtained in polynomial time.

We first check whether $\in X_0^*\bigcap \{x\mid~ l_i\le a_i^Tx\le u_i,~i=1,\ldots,m\}$
is empty, which can be done in polynomial time according to the following lemma.
\begin{lem}[\cite{HS}]
Let $A\in R^{m\times q}$ and $b\in R^m$, where $m$ is fixed and $q$
is arbitrary. For any given $r>0$, it is polynomially checkable
whether $\{u\in R^q\mid Au\le b, u^Tu=r\}$ is empty. Moreover, if
the set is nonempty, a feasible point can be found in polynomial
time.
\end{lem}
If $\in X_0^*\bigcap \{x\mid~ l_i\le a_i^Tx\le u_i,~i=1,\ldots,m\}\neq \emptyset$,
any point in $X_0^*\bigcap \{x\mid~ l_i\le a_i^Tx\le u_i,~i=1,\ldots,m\}$, globally solves $({\rm RS}_{m})$.
Otherwise,  we find  the unique local-nonglobal minimizer of $({\rm RS})$, denoted by $\overline{x}_0$, which is obtained in polynomial time according to Theorem \ref{thm:2}. Moreover, if
 \begin{equation}
 l_i< a_i^T\overline{x}_0< u_i,~i=1,2,\ldots,m, \label{interior}
 \end{equation}
then, $\overline{x}_0$ is the unique attained solution of the following problem:
\begin{eqnarray}
v({\rm RS}_{m}^0):= &\min & g(x)= \frac{1}{2}x^THx+c^Tx+\frac{\sigma}{4} \|x\|^4 \nonumber\\
&{\rm s.t.}~&l_i< a_i^Tx< u_i,~i=1,2,\ldots,m.\label{T_m_0}
\end{eqnarray}
It follows that
 \begin{equation}
v({\rm RS}_{m})=\min\{v({\rm RS}_{m}^0), v({\rm RS}_{m}^{1_1}), v({\rm RS}_{m}^{1_2}), \ldots, v({\rm RS}_{m}^{m_1}), v({\rm RS}_{m}^{m_2})\}
\label{iter1}
 \end{equation}
where $v({\rm RS}_{m}^{j_1})$ and $v({\rm RS}_{m}^{j_2})$ ($j=1,2,\ldots,m$) are defined as follows:
\begin{eqnarray}
v({\rm RS}_{m}^{j_1}):= &\min & g(x)= \frac{1}{2}x^THx+c^Tx+\frac{\sigma}{4} \|x\|^4 \nonumber\\
&{\rm s.t.}~& a_j^Tx= l_j, \label{y:e1}\\
&&l_i\le a_i^Tx\le u_i,~i=1,\ldots,j-1,j+1,\ldots,m, \nonumber
\end{eqnarray}
\begin{eqnarray}
v({\rm RS}_{m}^{j_2}):= &\min & g(x)= \frac{1}{2}x^THx+c^Tx+\frac{\sigma}{4} \|x\|^4 \nonumber\\
&{\rm s.t.}~& a_j^Tx= u_j, \label{y:e2}\\
&&l_i\le a_i^Tx\le u_i,~i=1,\ldots,j-1,j+1,\ldots,m. \nonumber
\end{eqnarray}
Otherwise, (\ref{interior}) does not hold true. Then, we have
\begin{equation}
v({\rm RS}_{m})=\min\{v({\rm RS}_{m}^{1_1}), v({\rm RS}_{m}^{1_2}), \ldots, v({\rm RS}_{m}^{m_1}), v({\rm RS}_{m}^{m_2})\}.
\label{iter2}
 \end{equation}
It remains to show how to solve $({\rm RS}_{m}^{j_1}),~j=1,\ldots,m$ as $({\rm RS}_{m}^{j_2})$ is similarly solved. Our
idea is to eliminate one variable using the equation (\ref{y:e1}) and
maintains the same structure as $({\rm RS}_{m})$.

Let $P_j\in \Bbb R^{n\times (n-1)}$ be a column-orthogonal matrix such
that $a_j^TP_j=0$. Let $z_0$ be a feasible solution to (\ref{y:e1}).
Then $z_0-P_jP_j^Tz_0$ is also feasible to (\ref{y:e1}). Using the
null-space representation, we have
\begin{equation}
\{x\in \Bbb R^n\mid a_j^Tx= b_j\}=\{z_0-P_jP_j^Tz_0+P_jz\mid z\in\Bbb  R^{n-1}\} \label{redu}
\end{equation}
and
\begin{eqnarray*}
\|x\|^4&=&\left((z_0-P_jP_j^Tz_0+P_jz)^T(z_0-P_jP_j^Tz_0+P_jz)\right)^2\\
&=&\left(z_0^T(I-P_jP_j^T)(I-P_jP_j^T)z_0+2z_0^T(I-P_jP_j^T)P_jz+z^TP_j^TP_jz\right)^2 \\
&=&\left(z_0^T(I-P_jP_j^T)z_0+z^Tz\right)^2.\\
&=&\left(z_0^T(I-P_jP_j^T)z_0\right)^2+2(z_0^T(I-P_jP_j^T)z_0)z^Tz+\|z\|^4
\end{eqnarray*}
We can
equivalently express $({\rm RS}_{m}^{j_1})$ as:
\begin{eqnarray*}
v({\rm RS}_{m}^{j_1}) = &\min & g(z_0-P_jP_j^Tz_0+P_jz)
%=\frac{1}{2}(z_0-P_jP_j^Tz_0+P_jz)^TH(z_0-P_jP_j^Tz_0+P_jz)+c^T(z_0-P_jP_j^Tz_0+P_jz)+
%\frac{\sigma}{4} (z_0^T(I-P_jP_j^T)z_0+z^Tz)^2
\\
&{\rm s.t.}&l_i\le a_i^T(z_0-P_jP_j^Tz_0+P_jz)\le u_i,~i=1,\ldots,m,~i\neq j
\end{eqnarray*}
which is again a special case of  $({\rm RS}_{m-1})$.

Iteratively applying (\ref{iter1}) or (\ref{iter2}), we will eventually terminate when
no linear constraint left. Let $s$ be the smallest
number such that any $s+1$ columns of $\{a_1,\ldots,a_m\}$ are dependent. By this
inductive way, there are at most $m\times(m-1)\times\cdots\times
(m-s+1)$  regularised subproblems to be solved.  Since $m$
is assumed to be fixed, the total number of reduction iterations is
bounded by a constant factor of $m$. We thus have proved that

\begin{thm}
For each fixed $m$, $({\rm RS}_m)$ is polynomially solvable.
\end{thm}

\section{Conclusions}

In this paper we have characterized the local and global minimizers of the regularised subproblem (RS) in optimization. We first show the existing necessary optimality condition  for (RS) in literature is also {\it sufficient} and the  $\ell_2$ norm of  the global minimizer is {\it always} unique. The hidden convexity of (RS) is also obtained. Then we establish a necessary and sufficient condition for the local-nonglobal minimizer of (RS). We notice that this condition remains open for the trust-region subproblem. As a corollary, we show (RS) with $p>2$ has at most {\it one} local-nonglobal minimizer. As a further application, we show (RS) with $p=4$ and a fixed number of linear inequality constraints can be solved in polynomial time, while general linear constrained (RS) is shown to be NP-hard. It is unknown what happens when $p\neq 4$.

%\section*{Acknowledgments}

\end{document}